\documentclass[12pt]{article}
\usepackage{amsfonts}
\usepackage{full page,
amssymb, amscd, amsmath,graphicx}

 \title{{\bf Vertex algebras associated to abelian current algebras}}
 \author{Jinwei Yang}
    \date{}
    \begin{document}
    \bibliographystyle{alpha}
    \maketitle
\newtheorem{thm}{Theorem}[section]
\newtheorem{defn}[thm]{Definition}
\newtheorem{prop}[thm]{Proposition}
\newtheorem{cor}[thm]{Corollary}
\newtheorem{lemma}[thm]{Lemma}
\newtheorem{rema}[thm]{Remark}
\newtheorem{app}[thm]{Application}
\newtheorem{prob}[thm]{Problem}
\newtheorem{conv}[thm]{Convention}
\newtheorem{conj}[thm]{Conjecture}
\newtheorem{cond}[thm]{Condition}
    \newtheorem{exam}[thm]{Example}
\newtheorem{assum}[thm]{Assumption}
     \newtheorem{nota}[thm]{Notation}
\newcommand{\halmos}{\rule{1ex}{1.4ex}}
\newcommand{\pfbox}{\hspace*{\fill}\mbox{$\halmos$}}
\newcommand{\nn}{\nonumber \\}

 \newcommand{\res}{\mbox{\rm Res}}
 \newcommand{\ord}{\mbox{\rm ord}}
\renewcommand{\hom}{\mbox{\rm Hom}}
\newcommand{\edo}{\mbox{\rm End}\ }
 \newcommand{\pf}{{\it Proof.}\hspace{2ex}}
 \newcommand{\epf}{\hspace*{\fill}\mbox{$\halmos$}}
 \newcommand{\epfv}{\hspace*{\fill}\mbox{$\halmos$}\vspace{1em}}
 \newcommand{\epfe}{\hspace{2em}\halmos}
  \newcommand{\nno}{\nonumber}
\newcommand{\nord}{\mbox{\scriptsize ${\circ\atop\circ}$}}
\newcommand{\wt}{\mbox{\rm wt}\ }
\newcommand{\swt}{\mbox{\rm {\scriptsize wt}}\ }
\newcommand{\lwt}{\mbox{\rm wt}^{L}\;}
\newcommand{\rwt}{\mbox{\rm wt}^{R}\;}
\newcommand{\slwt}{\mbox{\rm {\scriptsize wt}}^{L}\,}
\newcommand{\srwt}{\mbox{\rm {\scriptsize wt}}^{R}\,}
\newcommand{\clr}{\mbox{\rm clr}\ }
\newcommand{\tr}{\mbox{\rm Tr}}
\newcommand{\C}{\mathbb{C}}
\newcommand{\Z}{\mathbb{Z}}
\newcommand{\R}{\mathbb{R}}
\newcommand{\Q}{\mathbb{Q}}
\newcommand{\N}{\mathbb{N}}
\newcommand{\CN}{\mathcal{N}}
\newcommand{\F}{\mathcal{F}}
\newcommand{\I}{\mathcal{I}}
\newcommand{\V}{\mathcal{V}}
\newcommand{\one}{\mathbf{1}}
\newcommand{\BY}{\mathbb{Y}}
\newcommand{\ds}{\displaystyle}

        \newcommand{\ba}{\begin{array}}
        \newcommand{\ea}{\end{array}}
        \newcommand{\be}{\begin{equation}}
        \newcommand{\ee}{\end{equation}}
        \newcommand{\bea}{\begin{eqnarray}}
        \newcommand{\eea}{\end{eqnarray}}
         \newcommand{\lbar}{\bigg\vert}
        \newcommand{\p}{\partial}
        \newcommand{\dps}{\displaystyle}
        \newcommand{\bra}{\langle}
        \newcommand{\ket}{\rangle}

        \newcommand{\ob}{{\rm ob}\,}
        \renewcommand{\hom}{{\rm Hom}}

\newcommand{\A}{\mathcal{A}}
\newcommand{\Y}{\mathcal{Y}}

\renewcommand{\theequation}{\thesection.\arabic{equation}}
\renewcommand{\thethm}{\thesection.\arabic{thm}}
\setcounter{equation}{0} \setcounter{thm}{0}
\date{}
\maketitle

\begin{abstract}
We construct a family of vertex algebras associated to the current algebra of finite-dimensional abelian Lie algebras along with their modules and logarithmic modules. We show this family of vertex algebras and their modules are quasi-conformal and strongly $\N$-graded and verify convergence and extension property needed in the logarithmic tensor category theory for strongly graded logarithmic modules developed by Huang, Lepowsky and Zhang.
\end{abstract}

\section{Introduction}
In a series of papers \cite{HLZ1}-\cite{HLZ8}, Huang, Lepowsky and Zhang developed the theory of logarithmic tensor categories from strongly graded vertex algebras and their logarithmic modules. Representation theory of strongly graded vertex algebras has then been further studied in \cite{Y1} and \cite{Y2}. But so far, the only source of strongly graded vertex algebras and their modules comes from vertex algebras and modules associated with not necessarily positive definite even lattice. In this paper, we construct a new family of strongly graded vertex algebras along with a natural logarithmic module category and derive convergence and extension properties needed in the theory of logarithmic tensor categories.

The polynomial current algebras associated to finite simple Lie algebras have been studied in \cite{FL}, \cite{CG1} and \cite{CG2} et al. Current Lie algebra is the standard parabolic subalgebra of an affine Lie algebra and its representation has broad applications. There are numerous interesting and related families of modules in the representations of current Lie algebra such as Kirillov-Reshetikhin modules (see a detailed discussion in \cite{K}), Demazure modules arising from the positive level representations of the affine Lie algebra (see the details in \cite{BCM}), Weyl modules (\cite{CP}) and so on.

In this paper, we mainly focus on the current algebra of finite abelian Lie algebras, which is an infinite dimensional $\N$-graded abelian Lie algebra. Following the construction of vertex operator algebras associated to the Heisenberg Lie algebras (\cite{FLM}, \cite{LL}), we construct a family of vertex algebras $V$ associated to the affinization of the abelian current Lie algebras. This new family of vertex algebras does not possess a conformal vector and does not satisfy the grading restriction property required in the notion of vertex operator algebra. Instead, it has a second grading and satisfies the strongly gradedness condition required in the notion of strongly graded vertex algebra. Furthermore, we construct a series of operators which act on $V$ as generators $L(n)$ for $n \geq -1$ of the Virasoro algebra and hence the vertex algebra $V$ is a quasi-conformal vertex algebra in the sense of \cite{FB}.

We construct the strongly graded $V$-module category with irreducible and indecomposable objects associated to the evaluation modules for the current Lie algebras. Logarithmic intertwining operators among a triple of strongly graded logarithmic modules associated to the evaluation modules at $0$ can be constructed explicitly using the method developed in \cite{M2}. In general, it is still not clear how to construct logarithmic intertwining operators among logarithmic modules associated to the other families of modules for current Lie algebras.

Matrix elements of products and iterates of logarithmic intertwining operators among strongly graded logarithmic $V$-modules satisfy certain systems of differential equations under $C_1$-cofiniteness condition introduced in \cite{Y2} on the logarithmic $V$-modules. Using these systems of differential equations, we verify the convergence and extension property needed in the logarithmic tensor category theory for such strongly graded logarithmic $V$-modules.

This new family of vertex algebras provides very concrete examples for Huang-Lepowsky-Zhang's logarithmic tensor category theory. We will study the other conditions needed in the logarithmic tensor category, such as closedness under the $P(z)$-tensor product and associativity isomorphism, and build logarithmic tensor category for this new family of vertex algebras in the forthcoming papers.

This paper is structured as follows: In section 2, we introduce current algebra of a Lie algebra and its finite dimensional module category. In section 3 and 4, we construct a family of vertex algebras $V$ and their modules categories associated to current algebra of an abelian Lie algebra and its module categories. In section 5, we build logarithmic modules and construct logarithmic intertwining operators among triples of certain $V$-modules explicitly. In section 6, we study the strongly gradeness for $V$ and its module categories. In section 7, we derive differential equations for matrix elements of logarithmic intertwining operators and verify the convergence and expansion property needed in Huang-Lepowsky-Zhang's logarithmic tensor categories theory. In section 8, graded dimension for $V$ is calculated.

\paragraph{Acknowledgments}
I am grateful to Profs. Katrina Barron, Vyjayanthi Chari, Yi-Zhi Huang, Haisheng Li, Antun Milas and Ingo Runkel for valuable discussions.

\setcounter{equation}{0}
\section{Current algebra of a Lie algebra}
Current algebra of a finite dimensional Lie algebra has been studied in \cite{FL}, \cite{CG1}, \cite{CG2} et al. We first recall the definitions and notions, see for instance \cite{CG2}.

Let $\C[t]$ be the ring of polynomials in an indeterminate $t$. The current algebra $\mathfrak{g}[t]$ of a Lie algebra $\mathfrak{g}$ is the Lie algebra $\mathfrak{g}\otimes \C[t]$, where the Lie bracket is defined by
\[
[x \otimes f, y \otimes g]_{\mathfrak{g}[t]} = [x, y]_{\mathfrak{g}} \otimes fg,\;\;\; x, y \in \mathfrak{g},\;\; f, g \in \C[t].
\]
We will write $xf$ for the element $x\otimes f$, $x \in \mathfrak{g}$, $f \in \C[t]$ of $\mathfrak{g}[t]$.

In this paper, we will focus on the current algebra $\mathfrak{h}[t]$ of an finite dimensional abelian Lie algebra $\mathfrak{h}$ equipped with a nondegenerate symmetric bilinear form $\langle \cdot, \cdot \rangle_{\mathfrak{h}}$. The current Lie algebra $\mathfrak{h}[t]$ has trivial Lie bracket and has an invariant symmetric bilinear form induced from $\langle \cdot, \cdot \rangle_{\mathfrak{h}}$:
\[
\langle xt^m, yt^n\rangle_{\mathfrak{h}[t]} = \delta_{m,n}\langle x, y\rangle_{\mathfrak{h}}, \;\;\; x, y \in \mathfrak{h}, \;\; m, n \in \N.
\]
Let $\{u^{(1)}, \dots, u^{(d)}\}$ be an orthonormal basis of $\mathfrak{h}$ with respect to $\langle \cdot, \cdot \rangle_{\mathfrak{h}}$. Then $\{u^{(i)}t^j|1\leq i \leq d, j \in \N\}$ is an orthonormal basis of $\mathfrak{h}[t]$ with respect to $\langle \cdot, \cdot \rangle_{\mathfrak{h}[t]}$.

The form $\langle \cdot, \cdot \rangle_{\mathfrak{h}}$ being nondegenerate, we will identify $\mathfrak{h}$ with its dual space $\mathfrak{h}^{*}$. Let $\lambda \in \mathfrak{h} = \mathfrak{h}^{*}$. Denote by $\C_{\lambda}$ the one-dimensional $\mathfrak{h}$-module with $h \in \mathfrak{h}$ acting as the scalar $\langle\lambda, h\rangle$. For every $c \in \C$ we define an $\mathfrak{h}[t]$-module $V(\lambda, c) = \C_{\lambda}$ as a vector space  with action given by
\begin{equation}\label{ev module}
(hf)\cdot v = f(c)\lambda(h)v, \;\;\; h \in \mathfrak{h},\; f \in \C[t],\; v \in \C_{\lambda}.
\end{equation}

Finite-dimensional irreducible modules for $\mathfrak{h}[t]$ are well-known, see for instance \cite{FL}: Let $W$ be an irreducible module for $\mathfrak{h}[t]$. Then either $W$ is trivial or there exists $k \in \N$, $\lambda_i \in \mathfrak{h}^{*}$ and distinct $c_i \in \C$ for $1 \leq i \leq k$ such that
\begin{equation}\label{class ev module}
W = V(\lambda_1, c_1)\otimes \cdots \otimes V(\lambda_k, c_k).
\end{equation}

We define {\it generalized Casimir operator} $\Omega$ on $V(\lambda, c)$ for $|c|<1$ as follows:
\begin{equation}\label{casimir}
\Omega = \sum_{i=1}^d\sum_{n \in \N}(u^{(i)}t^n)(u^{(i)}t^n).
\end{equation}
Although there might be infinitely many nonzero terms in $\Omega\cdot w$ for $w \in V(\lambda, c)$, the sum is convergent for $|c| < 1$. In fact, $\Omega$ acts as the scalar $\frac{1}{1-c^2}\langle \lambda, \lambda \rangle_{\mathfrak{h}}$.

\setcounter{equation}{0}
\section{Vertex algebra associated to $\mathfrak{h}[t]$}
We will construct a vertex algebra structure associated to $\mathfrak{h}[t]$, following the steps in section $6.2$ and $6.3$ in \cite{LL}. The new vertex algebra generalizes the vertex operator algebra structure associated with Heisenberg Lie algebra.

In this paper, we will let $z, z_0, z_1, z_2$ denote commuting formal variables or complex variables.

Set
\[
\widehat{\mathfrak{h}[t]} = \mathfrak{h}[t]\otimes \C[s, s^{-1}] \oplus \C{\bf k},
\]
equipped with the bracket relations
\[
[xt^i\otimes s^m, yt^j\otimes s^n] = m\langle x, y\rangle_{\mathfrak{h}}\delta_{m+n,0}\delta_{i,j}{\bf k}
\]
for $x, y \in \mathfrak{h}$, $i, j \in \N$ and $m, n \in \Z$, together with the condition that ${\bf k}$ is a nonzero central element of $\widehat{\mathfrak{h}[t]}$.

The affine Lie algebra $\widehat{\mathfrak{h}[t]}$ is a $\Z$-graded Lie algebra with
\[
\widehat{\mathfrak{h}[t]} = \coprod_{n \in \Z}\widehat{\mathfrak{h}[t]}_{(n)},
\]
where
\[
\widehat{\mathfrak{h}[t]}_{(0)} = \mathfrak{h}[t]\oplus \C {\bf k} \;\;\;\mbox{and}\;\;\; \widehat{\mathfrak{h}[t]}_{(n)} = \mathfrak{h}[t]\otimes s^{-n}\;\;\; \mbox{for}\;\;\; n \neq 0.
\]
It has graded subalgebras
\[
\widehat{\mathfrak{h}[t]}_{+} = \mathfrak{h}[t]\otimes s^{-1}\C[s^{-1}],\;\;\;\mbox{and}\;\;\;
\widehat{\mathfrak{h}[t]}_{-} = \mathfrak{h}[t]\otimes s\C[s].
\]

Form the induced module
\[
M(l) = \mbox{Ind}_{\widehat{\mathfrak{h}[t]}_{-}\oplus \widehat{\mathfrak{h}[t]}_{(0)}}^{\widehat{\mathfrak{h}[t]}}(\C{\bf 1}) = U(\widehat{\mathfrak{h}[t]})\otimes_{U(\widehat{\mathfrak{h}[t]}_{-}\oplus \widehat{\mathfrak{h}[t]}_{(0)})} \C{\bf 1} = S(\widehat{\mathfrak{h}[t]}_{+})\otimes {\C}{\bf 1},
\]
where $\widehat{\mathfrak{h}[t]}_{-} \oplus \mathfrak{h}[t]$ annihilates ${\bf 1}$ and ${\bf k}$ acts as a scalar multiplication by $l$. Then $M(l)$ has a natural vertex algebra structure (see Theorem $6.2.11$ in \cite{LL}).

For $v = a_1(-n_1)\cdots a_k(-n_k){\bf 1} \in M(l)$ ($a_i \in \mathfrak{h}[t], n_i \in \Z_{+}$), the vertex operator map is given by
\[
Y(v, z) = \nord \frac{1}{(n_1-1)!}\big(\frac{d}{dz}\big)^{n_1-1}a_1(z)\cdots  \frac{1}{(n_k-1)!}\big(\frac{d}{dz}\big)^{n_k-1}a_k(z)\nord,
\]
where for $a \in \mathfrak{h}[t]$,
\[
a(z) = \sum_{n \in \Z}a(n)z^{-n-1}  \in (\edo M(l))[[z, z^{-1}]].
\]

We define operators $L(n)$ for $n \geq -1$ on $M(l)$ by
\begin{eqnarray}\label{def:Ln}
L(n) = \frac{1}{2l}\sum_{i=1}^d\sum_{j \in \N}\sum_{m \in \Z}\nord(u^{(i)}t^j)(n-m)(u^{(i)}t^j)(m)\nord.
\end{eqnarray}
The operators $L(n)$ are well-defined since for each $v \in M(l)$, $L(n)v$ has only finitely many nonzero terms.

We will show that the operators $L(n)$ for $n \geq -1$ have the following properties:

\begin{prop}
For $a \in \mathfrak{h}[t]$ and $k, n \in \Z, n \geq -1$,
\begin{equation}\label{e1}
[L(n), a(k)] = -ka(n+k)
\end{equation}
on $M(l)$. Furthermore, for $v = a_1(-n_1)\cdots a_k(-n_k){\bf 1} \in M(l)$ ($a_i \in \mathfrak{h}[t], n_i \in \Z_{+}$),
\begin{equation}\label{e2}
L(0)v = \big(\sum_{i=1}^k n_i\big)v
\end{equation}
\[
L(-1) = \mathcal{D},
\]
where $\mathcal{D}$ is the $\mathcal{D}$-operator of the vertex algebra $M(l)$.
\end{prop}
\pf We prove (\ref{e1}) for $a = ht^j$ ($h \in \mathfrak{h}, j \in \N$) and $k < 0$. The case $k \geq 0$ is similar. If $n \neq -2k$,
\begin{eqnarray*}
[L(n), (ht^j)(k)] &=& [\frac{1}{2l}\sum_{i=1}^d\sum_{j,m \in \N}(u^{(i)}t^j)(n-m)(u^{(i)}t^j)(m), (ht^j)(k)]\nn
&=& [\frac{1}{l}\sum_{i=1}^d(u^{(i)}t^j)(n+k)(u^{(i)}t^j)(-k), (ht^j)(k)]\nn
&=& -k\sum_{i=1}^d \langle u^{(i)}, h\rangle (u^{(i)}t^j)(n+k)\nn
&=& -k(ht^j)(n+k);
\end{eqnarray*}
 If $n = -2k$,
\begin{eqnarray*}
[L(n), (ht^j)(k)] &=& [\frac{1}{2l}\sum_{i=1}^d\sum_{j,m \in \N}(u^{(i)}t^j)(n-m)(u^{(i)}t^j)(m), (ht^j)(k)]\nn
&=& [\frac{1}{2l}\sum_{i=1}^d(u^{(i)}t^j)(-k)(u^{(i)}t^j)(-k), (ht^j)(k)]\nn
&=& -k\sum_{i=1}^d \langle u^{(i)}, h\rangle (u^{(i)}t^j)(-k)\nn
&=& -k(ht^j)(-k).
\end{eqnarray*}
Using (\ref{e1}) and the action of generalized Casimir operator on $M(l)$, we have
\begin{eqnarray*}
L(0)v &=& L(0)(a_1(-n_1)\cdots a_k(-n_k)v_{\lambda})\nn
&=& \big(\sum_{i=1}^kn_i\big) v + a_1(-n_1)\cdots a_k(-n_k)L(0){\bf 1}\nn
&=& \big(\sum_{i=1}^k n_i\big)v.
\end{eqnarray*}
From (\ref{e1}) and properties of operators $\mathcal{D}$, we have
\[
[L(-1)-\mathcal{D}, xt^j(n)] = 0
\]
as operators on $M(l)$, and we also have
\[
(L(-1) - \mathcal{D}){\bf 1} = 0.
\]
Since $M(l)$ is generated from ${\bf 1}$ by $U(\widehat{\mathfrak{h}[t]})$, it follows that $L(-1) = \mathcal{D}$ on $M(1)$. \epfv

\begin{thm}\label{quasiconformalproperty}
For $A \in M(l)$, $k \in \Z$ and $n \geq -1$, we have
\[
[L(n), Y(A, w)] = \sum_{m \geq -1}^n \binom{n+1}{m+1}w^{n-m}Y(L(m)A, w).
\]
In particular,
\[
[L(n), A_k] = \sum_{m = -1}^n\binom{n+1}{m+1}(L(m)A)_{k+n-m}.
\]
\end{thm}
\pf It suffices to prove the formula for homogeneous elements $A \in M(l)$. Let $A = h_1t^{i_1}(-n_1)\cdots h_kt^{i_k}(-n_k){\bf 1}$ ($h_i \in \mathfrak{h}, n_i \in \Z_{+}$). Set
\[
\omega = \sum_{i=1}^d\sum_{j = 1}^k(u^{(i)}t^{i_j})(-1)(u^{(i)}t^{i_j})(-1){\bf 1}
\]
and
\[
L'(n) = \res_z z^{n+1}Y(\omega, z).
\]
It is obvious that
\[
[L'(n), Y(A,w)] = [L(n), Y(A,w)], \;\;\; L'(n)A = L(n)A
\]
since if $i \neq j$, for $h_1, h_2 \in \mathfrak{h}$ and $m, n \in \Z$,
\[
[(h_1t^i)(m), (h_2t^j)(n)] = 0.
\]
From Theorem 8.6.1 in \cite{FLM}, we have
\begin{eqnarray}
[Y(\omega, z), Y(A, w)] &=& \res_{z_0}w^{-1}Y(Y(\omega, z_0)A, w)e^{-z_0(\partial/\partial z)}\delta(z/w)\nn
&=& \sum_{m \in \N}w^{-1}Y(L'(m-1)A, w)\frac{(-1)^m}{m!}(\frac{\partial}{\partial z})^m\delta(\frac{z}{w}).
\end{eqnarray}
By taking the coefficient of $z^{-n-2}$, we obtain for $n \geq -1$,
\begin{eqnarray*}
[L'(n), Y(A, w)] &=& \sum_{m = 0}^{n+1}\binom{n+1}{m}Y(L'(m-1)A, w)w^{n-m+1}\nn
&=& \sum_{m = -1}^{n}\binom{n+1}{m+1}Y(L'(m)A, w)w^{n-m}.
\end{eqnarray*}
Thus
\begin{equation*}
[L(n), Y(A, w)] = \sum_{m \geq -1}^n \binom{n+1}{m+1}w^{n-m}Y(L(m)A, w).
\end{equation*}\epfv

\begin{thm}\label{Virasoro relation}
For $n \geq -1$, the operators $L(n)$ satisfy the Virasoro Lie algebra relations. That is, for $m,n \geq -1$,
\[
[L(m), L(n)] = (m-n)L(m+n).
\]
\end{thm}
\pf For $n \geq -1$, set
\[
L_j(n) = \frac{1}{2l}\sum_{i=1}^d\sum_{m \in \Z}\nord(u^{(i)}t^j)(n-m)(u^{(i)}t^j)(m)\nord.
\]
Then
\begin{eqnarray*}
 [L(m), L(n)] &= & [\sum_{j\in \N}L_j(m), \sum_{j\in \N}L_j(n)] = \sum_{j \in \N}[L_j(m), L_j(n)] \nn
&= & (m-n)\sum_{j \in \N}L_j(m+n) = (m-n)L(m+n).
\end{eqnarray*}\epfv

We recall the following definition from \cite{FB}:
\begin{defn}{\rm
A vertex algebra $V$ is called {\it quais-conformal} if it carries the operators $L(n)$ for $n \geq -1$ such that for $m, n \geq -1$
\[
[L(m), L(n)] = (m-n)L(m+n)
\]
and for $v \in V$,
\[
[L(n), Y(v, x)] = \sum_{m \geq -1}^n \binom{n+1}{m+1}x^{n-m}Y(L(m)v, x).
\]}
\end{defn}

As an immediate consequence of Theorem \ref{quasiconformalproperty} and Theorem \ref{Virasoro relation}, we have:
\begin{cor}
The vertex algebra $M(l)$ is quasi-conformal.
\end{cor}

\setcounter{equation}{0}
\section{Restricted $\widehat{\mathfrak{h}[t]}$-modules}
In this section, we will construct certain modules for the quasi-conformal vertex algebra $M(l)$. First we introduce the following definitions.

We say that a $\widehat{\mathfrak{h}[t]}$-module $W$ is {\it restricted} $\Z$-graded of level $l$ if
\begin{enumerate}
\item[(i)]$W = \coprod_{n \in \Z}W_n$ such that
\[
a(n)W_m \subset W_{m-n}\;\;\; \mbox{for every}\; m,n \in \Z\; \mbox{and}\; a \in \mathfrak{h}[t];
\]
\item[(ii)]${\bf k}$ acts as multiplication by $l$ on $W$;
\item[(iii)]For every $w \in W$ and $n \in \Z_{+}$, $(xt^i)(n)w = 0$ for all but finitely many numbers of $i \in \N$;
\item[(iv)]There exists $N \in \N$ such that $W_n = 0$ for $n < N$.
\end{enumerate}
We denote by $\mathcal{C}_l$ the category of restricted $\Z$-graded $\widehat{\mathfrak{h}[t]}$-module of level $l$.

Let $V(\lambda, c)$ be the evaluation module at $c \in \C$ for $\mathfrak{h}[t]$ defined by (\ref{ev module}). Then the induced module
\[
W(\lambda, c, l) := M(l) \otimes V(\lambda, c)
\]
is a restricted $\Z$-graded $\widehat{\mathfrak{h}[t]}$-module of level $l$. In particular, the following proposition states that $W(\lambda, c, l)$ are simple objects in $\mathcal{C}_l$.

\begin{prop}\label{example}
Let $l$ be any nonzero complex number and let $\{u^{(1)}, \dots, u^{(d)}\}$ be an orthonormal basis of $\mathfrak{h}$. Let $\lambda \in \mathfrak{h} = \mathfrak{h}^{*}$. Set
\[
P(l, \lambda) = \C[x_{ijn}| i =1, \dots, d;\; j, n = 0, 1, 2, \dots],
\]
as a vector space, where $x_{ijn}$ are mutually commuting independent formal variables. Let $\widehat{\mathfrak{h}[t]}$ acts on $P(l, \lambda)$ by
\begin{eqnarray*}
{\bf k} &\longmapsto & l,\\
(u^{(i)}t^j)(0) &\longmapsto & c^{j}\langle u^{(i)}, \lambda\rangle,\\
(u^{(i)}t^j)(n) &\longmapsto & nl\frac{\partial}{\partial x_{ijn}},\\
(u^{(i)}t^j)(-n) &\longmapsto & x_{ijn}\;(\mbox{the multiplication operator})
\end{eqnarray*}
for $i = 1, \dots, d$ and $j, n \in \N$. Then $P(l, \lambda)$ is an irreducible $\widehat{\mathfrak{h}[t]}$-module. In particular, $P(l, \lambda) \cong W(\lambda, c, l)$ as irreducible $\widehat{\mathfrak{h}[t]}$-modules.
\end{prop}

In fact, all the restricted modules are of the following form:
\begin{thm}\label{main restricted}
Let $W$ be a restricted $\Z$-graded $\widehat{\mathfrak{h}[t]}$-module of level $l$. Then
\begin{equation}
W \cong M(l) \otimes \Omega(W),
\end{equation}
where $\Omega(W) = \{w \in W|\;a(n)w = 0, a \in \mathfrak{h}[t], n \in \Z_{+}\}$ (the vacuum space of $W$) is $\mathfrak{h}[t]$-stable.
\end{thm}
\pf We will prove the linear map
\begin{eqnarray*}
f: M(l) \otimes \Omega(W) &\longrightarrow & W\\
u \otimes w &\longmapsto & u\cdot w
\end{eqnarray*}
defines a $\widehat{\mathfrak{h}[t]}$-module isomorphism. First $f$ is injective. In fact, let $K$ be the kernel of $f$. Then $K \in \mathcal{C}_l$ and if $K \neq 0$ it contains a vacuum vector $w$. Then $w \in \Omega(W)$ since $\Omega(W)$ is precisely the vacuum space of $M(l) \otimes \Omega(W)$. But this contradicts the injectivity of $f$ on $\Omega(W)$.

Now we show that $f$ is surjective. Suppose instead that $W/{\rm Im}f \neq 0$. Then $W/{\rm Im}f \in \mathcal{C}_l$ and so contains a vacuum vector $v$. Let $w$ be a representative of $v$ in $W$. Then $w \notin {\rm Im}f$, $a(n)w \in {\rm Im}f$ for all $a \in \widehat{\mathfrak{h}[t]}$, $n \in \Z_{+}$ and there exists $N \in \Z_{+}$ such that $a(m)w = 0$ for all $m > N$. We will find $t \in {\rm Im}f$ such that
\begin{equation}\label{e3}
a(n)\cdot t = a(n)\cdot w
\end{equation}
for all $a \in \widehat{\mathfrak{h}[t]}$ and $n \in \Z_{+}$, since $t-w$ would then be a vacuum vector in $W$ but not in $\Omega(W)$, a contradiction. If suffices to prove equation (\ref{e3}) for all $a$ of the form $u^{(i)}t^j$ for $i = 1, \dots ,d$ and $j \in \N$.

Choose a basis $\{\omega_{\gamma}\}_{\gamma \in \Gamma}$ ($\Gamma$ an index set) of $\Omega(W)$, and note that
\[
{\rm Im}f = \coprod_{\gamma \in \Gamma}M(l)\otimes \C\omega_{\gamma},
\]
by the injectivity of $f$. For each $i = 1, \dots, d$, $j \in \N, n \in \Z_{+}$ and $\gamma \in \Gamma$, let $s_{ijn}(\gamma)$ be the component of $u^{(i)}t^j(n)\cdot w$ in $M(l)\otimes \C\omega_{\gamma}$ with respect to the decomposition. Then for all $i_1, i_2 = 1, \dots, d$ and $j_1, j_2 \in \N, n_1, n_2 \in \Z_{+}$, we have $u^{(i_1)}t^{j_1}(n_1)u^{(i_2)}t^{j_2}(n_2)\cdot w = u^{(i_2)}t^{j_2}(n_2)u^{(i_1)}t^{j_1}(n_1)\cdot w$, so that for all $\gamma \in \Gamma$,
\[
u^{(i_1)}t^{j_1}(n_1)\cdot s_{i_2j_2n_2}(\gamma) = u^{(i_2)}t^{j_2}(n_2)\cdot s_{i_1j_1n_1}(\gamma).
\]
Since $W$ is restricted, there exist a finite set $J \subset \N$ and a finite subset $\Gamma_0 \subset \Gamma$ such that $s_{ijn}(\gamma) = 0$ unless $j \in J$, $\gamma \in \Gamma_0$ and $n \leq N$. If we can find $t_{\gamma} \in M(l)\otimes \C\omega_{\gamma}$ such that $u^{(i)}t^j(n)\cdot t_{\gamma} = s_{ijn}(\gamma)$ for $i = 1, \dots, d$, $j \in \N, n \in \Z_{+}$ and $\gamma \in \Gamma_0$, then we can take $t = \sum_{\gamma \in \Gamma_0}t_{\gamma}$ and we will be done.

Fixing $\gamma \in \Gamma_0$ and identify $M(l)\otimes \C\omega_{\gamma}$ with the polynomial algebra on the generators $x_{ijn}$. Then
\[
\frac{\partial}{\partial x_{i_1j_1n_1}}s_{i_2j_2n_2}(\gamma) = \frac{\partial}{\partial x_{i_2j_2n_2}}s_{i_1j_1n_1}(\gamma)\;\;\; \mbox{for}\;i_1, i_2 = 1, \dots, d, j_1, j_2, n_1, n_2 \in \N.
\]
Recall that $s_{ijn}(\gamma) = 0$ for $j \notin J$ and for $n > N$, we see that each $s_{ijn}(\gamma)$ lies in the polynomial algebra for finitely many generators $x_{ijn}$ for $j \in J$ and for $n \leq N$. Thus there exists $s$ in this algebra such that
\[
l\frac{\partial}{\partial x_{ijn}}s = s_{ijn}(\gamma)
\]
for $j \in J$ and $n \leq N$ and hence for all $i,j,n$. We may therefore take $t_{\gamma} = s$. \epfv

\begin{thm}Let $l$ be any complex number and let $W$ be any restricted $\Z$-graded $\widehat{\mathfrak{h}[t]}$-module of level $l$. Then there exists a unique $M(l)$-module structure on $W$ such that for $a \in \mathfrak{h}[t]$,
\[
Y_W(a, z) = a_W(z) = \sum_{n \in \Z}a(n)z^{-n-1} \in (\edo W)[[z, z^{-1}]].
\]
The vertex operator map for this module structure is given by
\begin{eqnarray*}
&& Y_W(a_1(-n_1)\cdots a_k(-n_k){\bf 1}, z) \nn
&=& \nord \frac{1}{(n_1-1)!}\big(\frac{d}{dz}\big)^{n_1-1}(a_1)_W(z)\cdots  \frac{1}{(n_k-1)!}\big(\frac{d}{dz}\big)^{n_k-1}(a_k)_W(z)\nord 1_{W}
\end{eqnarray*}
for $k \geq 0$, $a_i \in \mathfrak{h}[t], n_i \in \Z_{+}$.
In particular, the modules $W(\lambda, c, l)$ for $|c|<1$ are irreducible modules for the quasi-conformal vertex algebra $M(l)$, and the operators $L(n)$ for $n \geq -1$ on $W(\lambda, c, l)$:
\begin{eqnarray}\label{def:LnW}
L(n)_W = \frac{1}{2l}\sum_{i=1}^d\sum_{j \in \N}\sum_{m \in \Z}\nord(u^{(i)}t^j)(n-m)(u^{(i)}t^j)(m)\nord
\end{eqnarray}
are well-defined and satisfy
\begin{equation}\label{maina}
[L(m)_W, L(n)_W] = (m-n)L(m+n)_W,
\end{equation}
\begin{equation}\label{mainb}
[L(n)_W, Y_W(A, w)] = \sum_{i \geq -1}^n \binom{n+1}{i+1}w^{n-i}Y_W(L(i)A, w)
\end{equation}
for $A \in M(l)$ and $m , n \geq -1$.
\end{thm}
\pf The first part is standard from Theorem 6.2.12 in \cite{LL}; the proof for the second part is the same as the proof for Theorem \ref{quasiconformalproperty} and Theorem \ref{Virasoro relation} in the vertex algebra case. \epfv

%\begin{cor}
%Let $W$ be an irreducible restricted $\Z$-graded $\widehat{\mathfrak{h}[t]}$-module of level $l$. Then
%\[
%W \simeq M(l) \otimes V(\lambda, 0)
%\]
%for some $\lambda \in \mathfrak{h}^{*}$.
%\end{cor}
%\pf By Theorem \ref{main restricted}, $W \simeq M(l) \otimes \Omega(W)$. Since $W$ is irreducible, $\Omega(W)$ is an irreducible $\mathfrak{h}[t]$-module. By (\ref{class ev module}) (cf. Proposition 3.9 \cite{CG1}),
%\[
%\Omega(W) \simeq V(\lambda_1, a_1)\otimes \cdots \otimes V(\lambda_k, a_k)
%\]
%for $k \in \N$, $\lambda_i \in \mathfrak{h}^{*}$ and distinct $a_i \in \C$. Because $W$ is restricted, we have $k = 1$ and $a_1 = 0$. That is, $\Omega(W) \simeq V(\lambda, 0)$ for some $\lambda \in \mathfrak{h}^{*}$. \epfv

\section{Logarithmic modules and logarithmic intertwining operators}
In this section, we will construct logarithmic modules for the quasi-conformal vertex algebra $M(l)$. First we introduce notations and definitions from \cite{HLZ1} and \cite{M1}.

\begin{defn}{\rm Let $V$ be a quasi-conformal vertex algebra. We say that a weak $V$-module $W$ is {\it logarithmic} if it admits a direct sum decomposition into generalized $L(0)$-eigenspaces, i.e. admits a Jordan form with respect to the action of $L(0)$. We say that a logarithmic module $W$ is {\it genuine} if under the action of $L(0)$ it admits at least one Jordan block of size $2$ or more.}
\end{defn}

In this section and next section, we assume dim $\mathfrak{h} = 1$ with a fixed nonzero vector $h$. Let $W$ be a $\Z$-graded restricted $\widehat{\mathfrak{h}[t]}$-module. Then $W \simeq M(l) \otimes \Omega(W)$ by Theorem \ref{main restricted}. We first have the following Lemma:
\begin{lemma}
Let $W$ be a restricted $\Z$-graded $\widehat{\mathfrak{h}[t]}$-module. Suppose that $\sum_{j\in \N}(ht^j)(0)^2|_{\Omega(W)}$ admits a Jordan block of size at least $2$. Then $W$ is a genuine logarithmic $\widehat{\mathfrak{h}[t]}$-module.
\end{lemma}
\pf It follows from the fact that $W \simeq M(l) \otimes \Omega(W)$ and $L(0)|_{\Omega(W)} = \frac{1}{2l}\sum_{j\in \N}(ht^j)(0)^2$. \epfv

Let $\lambda \in \mathfrak{h}^{*}$. Denote by $\Omega_{\lambda}$ a finite dimensional $\mathfrak{h}$-module such that
\[
(h - \lambda(h))^n|_{\Omega_{\lambda}} = 0
\]
for some $n \in \Z_{+}$. We define an evaluation $\mathfrak{h}[t]$-module $\Omega(\lambda, c)\simeq \Omega_{\lambda}$ as a vector space with action given by
\[
(hf)\cdot v = f(c)h\cdot v,
\]
where $f \in \C[t]$ and $v \in \Omega_{\lambda}$.

Form the induced module
\[
G(\lambda, c, l) = M(l)\otimes \Omega(\lambda, c)
\]
for some $\lambda \in \mathfrak{h}^{*}$ and $|c| < 1$. We have that
\[
L(0)|_{\Omega(\lambda, c)} = \frac{1}{2l}\sum_{j \in \N}(ht^j)(0)^2 = \frac{1}{2l}\big(\sum_{j \in \N}c^{2j}\big)h^2 =  \frac{1}{2l(1-c^{2})}h^2.
\]
Then $G(\lambda, c, l)$ is a genuine logarithmic $\widehat{\mathfrak{h}[t]}$-module if $h(0)^2|_{\Omega(\lambda, c)}$ admits a Jordan block of size at least $2$.

\begin{defn}{\rm
Let $(W_1,Y_1)$, $(W_2,Y_2)$ and $(W_3,Y_3)$ be logarithmic modules
for a quasi-conformal vertex algebra $V$. A {\em logarithmic intertwining
operator of type ${W_3\choose W_1\,W_2}$} is a linear map
\begin{equation*}
{\cal Y}(\cdot, z)\cdot: W_1\otimes W_2\to W_3[\log z]\{z\},
\end{equation*}
or equivalently,
\begin{equation*}
w_{(1)}\otimes w_{(2)}\mapsto{\cal Y}(w_{(1)},z)w_{(2)}=\sum_{n\in
{\mathbb C}}\sum_{k\in {\mathbb N}}{w_{(1)}}_{n;\,k}^{\cal
Y}w_{(2)}z^{-n-1}(\log x)^k\in W_3[\log z]\{z\}
\end{equation*}
for all $w_{(1)}\in W_1$ and $w_{(2)}\in W_2$, such that the
following conditions are satisfied: the {\em lower truncation
condition}: for any $w_{(1)}\in W_1$, $w_{(2)}\in W_2$ and $n\in
{\mathbb C}$,
\begin{equation*}
{w_{(1)}}_{n+m;\,k}^{\cal Y}w_{(2)}=0\;\;\mbox{ for }\;m\in {\mathbb
N} \;\mbox{ sufficiently large,\, independently of}\;k;
\end{equation*}
the {\em Jacobi identity}:
\begin{eqnarray*}
\lefteqn{\dps z^{-1}_0\delta \bigg( {z_1-z_2\over z_0}\bigg)
Y_3(v,z_1){\cal Y}(w_{(1)},z_2)w_{(2)}}\nno\\
&&\hspace{2em}- z^{-1}_0\delta \bigg( {z_2-z_1\over -z_0}\bigg)
{\cal Y}(w_{(1)},z_2)Y_2(v,z_1)w_{(2)}\nno \\
&&{\dps = z^{-1}_2\delta \bigg( {z_1-z_0\over z_2}\bigg) {\cal
Y}(Y_1(v,z_0)w_{(1)},z_2) w_{(2)}}
\end{eqnarray*}
for $v\in V$, $w_{(1)}\in W_1$ and $w_{(2)}\in W_2$; the {\em $L(-1)$-derivative property:} for any
$w_{(1)}\in W_1$,
\begin{equation*}
{\cal Y}(L(-1)w_{(1)},z)=\frac d{dz}{\cal Y}(w_{(1)},z).
\end{equation*}}\epfv
\end{defn}

%\begin{defn}\label{slog:def}{\rm In the setting of Definition \ref{log:def}, suppose in addition that $V$ and $W_1$, $W_2$ and $W_3$ are strongly graded. A logarithmic intertwining operator $\cal{Y}$ as in Definition \ref{log:def} is a {\it grading-compatible logarithmic intertwining operator} if for $\beta, \gamma \in \tilde{A}$ and $w_1 \in W_1^{(\beta)}$, $w_2 \in W_2^{(\gamma)}$, $n \in \C$ and $k \in \N$, we have
%\[
%(w_1)_{n; k}w_2 \in W_3^{(\beta + \gamma)}.
%\]
%}
%\end{defn}

%\begin{defn}
%{\rm In the setting of Definition \ref{slog:def}, the grading-compatible logarithmic intertwining operators of a fixed type ${W_3\choose W_1\,W_2}$ form a vector space, which we denote by $\mathcal{V}_{W_1 W_2}^{W_3}$. We call the dimension of $\mathcal{V}_{W_1 W_2}^{W_3}$ the {\it fusion rule} for $W_1$, $W_2$ and $W_3$ and denote it by $N_{W_1 W_2}^{W_3}$.}
%\end{defn}

Following the construction in \cite{M2}, we can identify the space of logarithmic intertwining operators among triples of logarithmic modules $W_i = G(\lambda_i, 0, l)$ for $\lambda_i \in \mathfrak{h}^{*}$ ($i = 1, 2, 3$) with $\hom_{\mathfrak{h}[t]}(\Omega(\lambda_1, 0), \hom(\Omega(\lambda_2, 0), \Omega(\lambda_3, 0))$.

\setcounter{equation}{0}
\section{Strongly graded vertex algebras and their modules}
In this section, we define the notions of strongly $\N$-graded vertex algebra and its strongly $\N$-graded modules following \cite{HLZ1} (cf. \cite{Y1}, \cite{Y2}).
\begin{defn}\label{def:dgv}
{\rm A quasi-conformal vertex algebra
\[
V=\coprod_{n\in {\mathbb Z}} V_{(n)}
\]
is said to be {\em strongly graded with respect to $\N$} (or {\em
strongly $\N$-graded}) if it is equipped with a second grading, by $\N$,
\[
V=\coprod _{m \in \N} V^{(m)},
\]
such that the following conditions are satisfied: the two gradings
are compatible, that is,
\[
V^{(m)}=\coprod_{n\in {\mathbb Z}} V^{(m)}_{(n)}, \;\;
\mbox{where}\;V^{(m)}_{(n)}=V_{(n)}\cap V^{(m)}\;
\mbox{ for any }\;m \in \N;
\]
for any $m,k\in \N$ and $n\in {\mathbb Z}$,
\begin{eqnarray}
&&V^{(m)}_{(n)}=0\;\;\mbox{ for }\;n\;\mbox{ sufficiently
negative};\nn
&&\dim V^{(m)}_{(n)} <\infty;\nn
&&{\bf 1}\in V^{(0)}_{(0)};\nn
&& \label{Virasoro grading of v}L(j)V^{(m)} \subset V^{(m)}\;\;\;\;\mbox{for}\;\; j \geq -1;\\
\label{stronglygradingofv}&&v_j V^{(k)} \subset \coprod_{0\leq i \leq m+k}V^{(i)}\;\; \mbox{ for any
}\;v\in V^{(m)},\;j \in {\mathbb Z}.
\end{eqnarray}
}
\end{defn}

Note here the definition of the notion of strongly $\N$-graded vertex algebra is slightly different from \cite{HLZ1} where (\ref{stronglygradingofv}) is replaced by a stronger condition:
\[
v_l V^{(k)} \subset V^{(m+k)}\;\; \mbox{ for any
}\;v\in V^{(m)},\;l\in {\mathbb Z}.
\]

For modules for a strongly $\N$-graded algebra, we will also have a second
grading by $\N$.

\begin{defn}\label{def:dgw}{\rm
Let $V$ be a strongly $\N$-graded
quasi-conformal vertex algebra. A $V$-module
\[
W=\coprod_{n\in{\C}} W_{(n)}
\]
is said to be {\em strongly graded with respect to $\N$} (or
{\em strongly $\N$-graded}) if it is equipped with a
second gradation, by $\N$,
\[
W=\coprod _{m \in \N} W^{(m)},
\]
such that the following conditions are satisfied: the two gradations
are compatible, that is, for any $m \in \N$,
\[
W^{(m)}=\coprod_{n\in {\C}} W^{(m)}_{(n)},
\;\;\mbox{where }\; W^{(m)}_{(n)}=W_{(n)}\cap W^{(m)};
\]
for any $m, k \in \N$ and $n\in {\C}$,
\begin{eqnarray}
&&W^{(m)}_{(n)}=0 \;\;
\mbox{for}\; n\;\mbox{with
sufficiently
negative real part};\nn
&&\dim W^{(m)}_{(n)} < \infty; \nn
&& \label{Virasoro grading of w}L(j)W^{(m)} \subset W^{(m)}\;\;\;\;\mbox{for}\;\; j \geq -1;\\
\label{stronglygrading}&&v_j W^{(k)} \subset \coprod_{0\leq i \leq m+k}W^{(i)}\;\;\mbox{ for any
}\;v\in V^{(m)},\;j\in {\mathbb Z}.
\end{eqnarray}}
\end{defn}

Note here for the purpose of this paper, we modify the notion of strongly
$\N$-graded modules in \cite{HLZ1}, where (\ref{stronglygrading}) is replaced by
\[
v_l W^{(k)} \subset W^{(m+k)}\;\;\mbox{ for any
}\;v\in V^{(m)},\;l\in {\mathbb Z}.
\]

Similarly, we can define the notion of {\it strongly $\N$-graded logarithmic modules} by replacing $L(0)$-eigenspace $W_{(n)}$ by generalized $L(0)$-eigenspace $W_{[n]}$.

We define an $\N$-grading for $M(l)$ by
\[
\N\text{-}\wt x_1t^{i_1}(-n_1)\cdots x_kt^{i_k}(-n_k){\bf 1} = i_1 + \cdots + i_k.
\]
It is easy to check that $M(l)$ satisfies all the assumptions in Definition {\ref{def:dgv}}. Thus we have
\begin{cor}
The vertex algebra $M(l)$ is a strongly $\N$-graded quasi-conformal vertex algebra.
\end{cor}

We define an $\N$-grading for $W(\lambda, c, l)$ by
\[
\N\text{-}\wt x_1t^{i_1}(-n_1)\cdots x_kt^{i_k}(-n_k)v_{\lambda} = i_1 + \cdots + i_k,
\]
where $v_{\lambda}$ is a nontrivial element in $V(\lambda, c)$. Similarly, we define an $\N$-grading for $G(\lambda, c, l)$ by
\[
\N\text{-}\wt x_1t^{i_1}(-n_1)\cdots x_kt^{i_k}(-n_k)w = i_1 + \cdots + i_k,
\]
where $w \in \Omega(\lambda, c)$. Thus we construct strongly $\N$-graded modules and logarithmic modules for $M(l)$:
\begin{cor}For $\lambda \in \mathfrak{h}^{*}$ and $|c|<1$, let $V(\lambda, c)$ and $\Omega(\lambda, c)$ be $\mathfrak{h}[t]$-modules defined as before. Then $W(\lambda, c, l)$ and $G(\lambda, c, l)$ are strongly $\N$-graded modules and logarithmic modules for $M(l)$, respectively.
\end{cor}

In the rest of this section, we will recall the following useful definitions and notations from \cite{HLZ1} (cf. \cite{Y1}, \cite{Y2}):
\begin{defn}\label{DH}{\rm Let $V$ be a strongly $\N$-graded conformal vertex
algebra. The subspaces $V_{(n)}^{(m)}$ for $n \in \mathbb{Z}$, $m
\in \N$ are called the {\it doubly
homogeneous subspaces} of $V$. The elements in $V_{(n)}^{(m)}$
are called {\it doubly homogeneous} elements. Similar definitions
can be used for $W^{(m)}_{(n)}$ (respectively,  $W^{(m)}_{[n]}$) in the strongly graded (logarithmic) module $W$.}
\end{defn}

\begin{nota}\label{A-wt}
{\rm Let $v$ be a doubly homogeneous element of $V$. Let wt $v_n$, $n \in \mathbb{Z}$, refer to the weight of $v_n$ as an operator acting on $W$, and let $\N$-wt $v_n$ refer to the $\N$-weight of $v_n$ on $W$. Similarly, let $w$ be a doubly homogeneous element of $W$. We use wt $w$ to denote the weight of $w$ and $\N$-wt $w$ to denote the $\N$-grading of $w$.}
\end{nota}

\begin{lemma}\label{b}
Let $v \in V^{(m)}_{(n)}$, for $n \in \mathbb{Z}$,
$m \in \N$. Then for $k \in \mathbb{Z}$, {\rm wt}$\ v_k$ = $n - k -1$ and $\N$-{\rm wt}$\ v_k \leq m$.
\end{lemma}
\pf The first equation is standard from the theory of graded conformal vertex algebras
and the second follows easily from the definitions. \epfv

With the strong gradedness condition on a (logarithmic) module, we can now define the corresponding notion of contragredient (logarithmic) module.
\begin{defn}
{\rm
Let $W = \coprod_{m \in \N, n \in \C} W_{[n]}^{(m)}$ be a strongly $\N$-graded logarthmic
module for a strongly $\N$-graded quasi-conformal vertex algebra. For each $m \in \N$ and $n \in \C$, let us
identify $(W_{[n]}^{(m)})^{*}$ with the subspace of $W^{*}$ consisting of the linear function on $W$ vanishing
on each $W_{[n']}^{(m')}$ with $m' \neq m$ or $n' \neq n$. We define $W'$ to be the $(\N \times \C)$-graded
vector subspaces of $W^{*}$ given by
\[
W' =  \coprod_{m \in \N, n \in \C} (W')_{[n]}^{(m)}, \;\; \mbox{where}\;\;  (W')_{[n]}^{(m)} =  (W_{[n]}^{(-m)})^{*}.
\]
}
\end{defn}

The {\it adjoint vertex operators} $Y'(v, z)\; (v \in V)$ on $W'$ is defined in the same way as vertex operator algebra in section 5.2 in \cite{FHL}:
\begin{equation}\label{adjoint}
\langle Y'(v, z)w', w\rangle = \langle w', Y(e^{zL(1)}(-z^{-2})^{L(0)}v, z^{-1})w\rangle
\end{equation}
for $w' \in W', w\in W$. The pair $(W', Y')$ carries a strongly $-\N$-graded module structure.

For $u \in V$ and $n \in \Z$, let $u_n = \res_z z^nY'(v, z)$ be the operator on $W'$ and let $u_n^{*}: W \rightarrow W$ be the adjoint of $u_n$. From formula (\ref{adjoint}), we have
\[
\wt u_n^{*} = -\wt u_n
\]
and
\begin{equation}\label{N-wt}
\N\text{-}\wt u_n^{*} = \N\text{-}\wt u_n.
\end{equation}

\setcounter{equation}{0}
\section{$C_1$-cofiniteness condition and differential equations}
In this section, we let $V$ be the strongly $\N$-graded quasi-conformal vertex algebra $M(l)$ and we assume every strongly $\N$-graded $V$-module has weights in $\R$.

\begin{defn}{\rm
Let $W$ be a strongly $\N$-graded $V$-module and let $C_1(W)$ be the subspace of $W$ spanned by elements of the form $u_{-1}w$ for $u \in V_{+} = \coprod_{n>0}V_{(n)}$ and $w \in W$. The $\N$-grading on $W$ induces an $\N$-grading on $W/C_1(W)$:
\[
W/C_1(W) = \coprod_{m \in \N}(W/C_1(W))^{(m)},
\]
where
\[
(W/C_1(W))^{(m)} = W^{(m)}/(C_1(W))^{(m)}
\]
for $m \in \N$. If dim $W^{(m)}/(C_1(W))^{(m)} < \infty$ for $m \in \N$, we say that $W$ is {\it $C_1$-cofinite with respect to $\N$} or $W$ satisfies the {\it $C_1$-cofiniteness condition with respect to $\N$.}}
\end{defn}

\begin{cor}
For $\lambda \in \mathfrak{h}^{*}$ and $c \in \C$, let $V(\lambda, c)$ and $\Omega(\lambda, c)$ be $\mathfrak{h}[t]$-modules defined as before. Then $W(\lambda, c, l)$ and $G(\lambda, c, l)$ satisfy the $C_1$-cofiniteness condition with respect to $\N$.
\end{cor}

Let $R = \C[z_1^{\pm 1}, z_2^{\pm 1}, (z_1-z_2)^{-1}]$, $W_i, i=0,1,2,3$ be strongly $\N$-graded $V$-modules satisfying the $C_1$-cofiniteness condition with respect to $\N$ and
\[
T = R \otimes W_0 \otimes W_1 \otimes W_2 \otimes W_3
\]
which has a natural $R$-module structure. For simplicity, we shall omit one tensor symbol to write $f(z_1,z_2)\otimes w_0\otimes w_1\otimes w_2 \otimes w_3$ as $f(z_1,z_2)w_0\otimes w_1\otimes w_2 \otimes w_3$ in $T$.

For $u \in V_{+}$ and $w_i \in W_i$, $i = 0,1,2,3$, let $J$ be the submodule of $T$ generated by elements of the form
\begin{eqnarray*}
&&\mathcal{A}(u, w_0, w_1, w_2, w_3)\nno \\
&=& \sum_{k \geq 0}\left(\begin{array}{c}-1\\ k \end{array}\right)(-z_1)^ku_{-1-k}^{*}w_0 \otimes w_1 \otimes w_2 \otimes w_3 - w_0 \otimes u_{-1}w_1 \otimes w_2 \otimes w_3 \nno\\
&& -\sum_{k \geq 0}\left(\begin{array}{c}-1\\ k \end{array}\right)(-(z_1 - z_2))^{-1-k}w_0 \otimes w_1 \otimes u_kw_2 \otimes w_3\nno\\
&& -\sum_{k \geq 0}\left(\begin{array}{c}-1\\ k \end{array}\right)(-z_1)^{-1-k}w_0 \otimes w_1 \otimes w_2 \otimes u_kw_3,\nno\\
&&\mathcal{B}(u, w_0, w_1, w_2, w_3)\nno \\
&=& \sum_{k \geq 0}\left(\begin{array}{c}-1\\ k \end{array}\right)(-z_2)^ku_{-1-k}^{*}w_0 \otimes w_1 \otimes w_2 \otimes w_3 \nno\\
&& -\sum_{k \geq 0}\left(\begin{array}{c}-1\\ k \end{array}\right)(-(z_1 - z_2))^{-1-k}w_0 \otimes u_kw_1 \otimes w_2 \otimes w_3  - w_0 \otimes w_1 \otimes u_{-1}w_2 \otimes w_3\nno\\
&& -\sum_{k \geq 0}\left(\begin{array}{c}-1\\ k \end{array}\right)(-z_2)^{-1-k}w_0 \otimes w_1 \otimes w_2 \otimes u_kw_3,\nno\\
&&\mathcal{C}(u, w_0, w_1, w_2, w_3)\nno \\
&=& u_{-1}^{*}w_0 \otimes w_1 \otimes w_2 \otimes w_3 - \sum_{k \geq 0}\left(\begin{array}{c}-1\\ k \end{array}\right)z_1^{-1-k}w_0 \otimes u_kw_1 \otimes w_2 \otimes w_3 \nno\\
&& -\sum_{k \geq 0}\left(\begin{array}{c}-1\\ k \end{array}\right)z_2^{-1-k}w_0 \otimes w_1 \otimes u_kw_2 \otimes w_3 - w_0 \otimes w_1 \otimes w_2 \otimes u_{-1}w_3,\nno\\
&&\mathcal{D}(u, w_0, w_1, w_2, w_3)\nno \\
&=& u_{-1}w_0 \otimes w_1 \otimes w_2 \otimes w_3 \nno\\
&& -\sum_{k \geq 0}\left(\begin{array}{c}-1\\ k \end{array}\right)z_1^{k+1}w_0 \otimes e^{z_1^{-1}L(1)}(-z_1^2)^{L(0)}u_k(-z_1^{-2})^{L(0)}e^{-z_1^{-1}L(1)}w_1 \otimes w_2 \otimes w_3\nno\\
&& -\sum_{k \geq 0}\left(\begin{array}{c}-1\\ k \end{array}\right)z_2^{k+1}w_0 \otimes w_1 \otimes e^{z_2^{-1}L(1)}(-z_2^2)^{L(0)}u_k(-z_2^{-2})^{L(0)}e^{-z_2^{-1}L(1)}w_2 \otimes w_3\nno\\
&& -w_0 \otimes w_1 \otimes w_2 \otimes u_{-1}^{*}w_3.
\end{eqnarray*}

The double gradings on $W_i$ for $i = 0,1,2,3$ induce double gradings on $W_0\otimes W_1\otimes W_2 \otimes W_3$ and then also on $T$ (here we define the double gradings of elements of $R$ to be $0$). Let $T_{(r)}^{(m)}$ be the doubly homogeneous subspaces of weight $r$ and $\N$-weight $m$ for $r \in \R$ and $m \in \N$. Then
\[
T = \coprod_{r \in \R, m \in \N}T^{(m)}_{(r)},
\]
where $T^{(m)}_{(r)}$ are finitely generated $R$-modules and $T^{(m)}_{(r)} = 0$ when $r$ is sufficiently small. For $m \in \N, r \in \R$, we define the following filtration on $T$:
\begin{eqnarray*}
&& F_r(T) = \coprod_{s\leq r}T_{(s)}\nn
&& F^m(T) = \coprod_{0\leq k\leq m}T^{(k)}\nn
&& F_r^m(T) = \coprod_{s\leq r, 0\leq k\leq m}T^{(k)}_{(s)}.
\end{eqnarray*}
Then $F_r^m(T)$ are finitely generated $R$-modules, and
\[
F_r^m(T) \subset F_s^m(T)\;\; \mbox{for}\; r \leq s,
\]
\[
F^m(T) = \cup_{r \in \R}F_r^m(T).
\]
For $m \in \N, r \in \R$, we define the following filtration on $J$ induced from the above filtration on $T$:
\begin{eqnarray*}
&& F_r(J) = J \cap F_r(T)\nn
&& F^m(J) = J \cap F^m(T)\nn
&& F_r^m(J) = J \cap F_r^m(T).
\end{eqnarray*}
Then $F_r^m(J)$ are finitely generated $R$-modules, and
\[
F_r^m(J) \subset F_s^m(J)\;\;\; \mbox{for}\; r \leq s,
\]
\[
J^{(m)} = \cup_{r \in \R}F_r^m(J).
\]

\begin{prop}There exists $M \in \Z$ such that for any $r \in \R, m \in \N$, $F_r^m(T) \subset F_r^m(J) + F_M^m(T)$. In particular, $F^m(T) = F^m(J) + F_{M}^m(T)$.
\end{prop}
\pf Since dim $W_i^{(m_i)}/(C_1(W_i))^{(m_i)} < \infty$ for $m_i \in \N, i = 0, 1, 2, 3$, there exists $M \in \Z$ such that
\begin{eqnarray}\label{maineq}
\coprod_{n > M}(F^m(T))_{(n)} &\subset & \coprod_{m_0+m_1+m_2+m_3 \leq m, m_i \in \N}R((C_1(W_0))^{(m_0)}\otimes W_1^{(m_1)} \otimes W_2^{(m_2)}\otimes W_3^{(m_3)})\nn
&& +\; R(W_0^{(m_0)}\otimes (C_1(W_1))^{(m_1)} \otimes W_2^{(m_2)}\otimes W_3^{(m_3)})\nn
&& +\; R(W_0^{(m_0)}\otimes W_1^{(m_1)} \otimes (C_1(W_2))^{(m_2)}\otimes W_3^{(m_3)})\nn
&& +\; R(W_0^{(m_0)}\otimes W_1^{(m_1)} \otimes W_2^{(m_2)}\otimes (C_1(W_3))^{(m_3)}).
\end{eqnarray}
Note the right-hand side of (\ref{maineq}) is a finite sum.

We use induction on $r \in \R$. If $r$ is equal to $M$, $F_M^m(T) \subset F_M^m(J) + F_M^m(T)$. Now we assume that $F_r^m(T) \subset F_r^m(J) + F_M^m(T)$ for $r < s$ where $s>M$. We want to show that any doubly homogeneous elements of $T^{(m)}_{(s)}$ can be written as a sum of element of $F_s^m(J)$ and an element of $F_M^m(T)$. Since $s > M$, by (\ref{maineq}), any element of $T^{(m)}_{(s)}$ is an element of the right-hand side of (\ref{maineq}). We shall discuss only the case that this element is in $R((C_1(W_0))^{(m_0)}\otimes W_1^{(m_1)} \otimes W_2^{(m_2)}\otimes W_3^{(m_3)})$ for some $m_i \in \N$ such that $m_0+m_1+m_2+m_3 \leq m$; the other cases are completely similar.

We need only discuss elements of the form $w_0 \otimes u_{-1}w_1 \otimes w_2 \otimes w_3$, where $w_i \in W_i^{(m_i)}$ for $i = 0, 1, 2, 3$ and $u \in V_{+}$. By assumption, the weight of $w_0 \otimes u_{-1}w_1 \otimes w_2 \otimes w_3$ is $s$, then the weight of $u_{-1-k}^{*}w_0 \otimes w_1 \otimes w_2 \otimes w_3$, $w_0 \otimes w_1 \otimes u_kw_2 \otimes w_3$ and $w_0 \otimes w_1 \otimes w_2 \otimes u_kw_3$ for $k \geq 0$, are all less than $s$. Note that all these elements have $\N$-weights less than or equal to $m$ by (\ref{stronglygrading}) and (\ref{N-wt}). Thus these elements lie in $F_{s-1}^m(T)$. Also, since $\mathcal{A}(u, w_0, w_1, w_2, w_3) \in F_s^m(J)$, we see that
\begin{eqnarray*}
&&w_0 \otimes u_{-1}w_1 \otimes w_2 \otimes w_3\nno \\
&& \;\;\;\;\;\;\;\;= \mathcal{A}(u, w_0, w_1, w_2, w_3) + \sum_{k \geq 0}\left(\begin{array}{c}-1\\ k \end{array}\right)(-z_1)^ku_{-1-k}^{*}w_0 \otimes w_1 \otimes w_2 \otimes w_3 \nno\\
&& \;\;\;\;\;\;\;\;\;\;\;\; -\sum_{k \geq 0}\left(\begin{array}{c}-1\\ k \end{array}\right)(-(z_1 - z_2))^{-1-k}w_0 \otimes w_1 \otimes u_kw_2 \otimes w_3\nno\\
&& \;\;\;\;\;\;\;\;\;\;\;\; -\sum_{k \geq 0}\left(\begin{array}{c}-1\\ k \end{array}\right)(-z_1)^{-1-k}w_0 \otimes w_1 \otimes w_2 \otimes u_kw_3
\end{eqnarray*}
can be written as a sum of an element of $F_s^m(J)$ and elements of $F_{s-1}^m(T)$. Thus by the induction assumption, the element $w_0 \otimes u_{-1}w_1 \otimes w_2 \otimes w_3$ can be written as a sum of an element of $F_s^m(J)$ and an element of $F_M^m(T)$.

Now we have
\begin{eqnarray*}
F^m(T) &=& \coprod_{r \in \R}F_r^m(T)\nn
&\subset &\coprod_{r \in \R}F_r^m(J) + F_M^m(T)\nn
& =& F^m(J) + F_M^m(T).
\end{eqnarray*}
But we know that $F^m(J) + F_M^m(T) \subset F^m(T)$. Thus we have $F^m(T) = F^m(J) + F_M^m(T)$.
\epfv

We immediately have the following:
\begin{cor}\label{maincor}
For each $m \in \N$, the quotient $R$-module $F^m(T)/F^m(J)$ is finitely generated.
\end{cor}

For an element $\mathcal{W} \in F^m(T)$, we shall use $[\mathcal{W}]$ to denote the equivalence class in $F^m(T)/F^m(J)$ containing $\mathcal{W}$. We have the following theorem:

\begin{thm}
Let $W_i$ be strongly $\N$-graded generalized $V$-modules for $i = 0, 1, 2, 3$. For any $w_i \in W_i$ $(i = 0, 1, 2, 3)$, let $M_1$ and $M_2$ be the $R$-submodules of $F(T)/F(J)$ generated by $[w_0 \otimes L(-1)^jw_1 \otimes w_2 \otimes w_3]$, $j \geq 0$, and by $[w_0 \otimes w_1 \otimes L(-1)^jw_2 \otimes w_3]$, $j \geq 0$, respectively. Then $M_1$, $M_2$ are finitely generated. In particular, for any $w_i \in W_i$ $(i = 0, 1, 2, 3)$, there exist $a_k(z_1, z_2)$, $b_l(z_1, z_2) \in R$ for $k = 1, \dots, m$ and $l = 1, \dots, n$ such that
\begin{eqnarray}\label{exe1}
[w_0 \otimes L(-1)^mw_1 \otimes w_2 \otimes w_3] + a_1(z_1, z_2)[w_0 \otimes L(-1)^{m-1}w_1 \otimes w_2 \otimes w_3]\nn
+ \cdots + a_m(z_1, z_2)[w_0 \otimes w_1 \otimes w_2 \otimes w_3] = 0,
\end{eqnarray}

\begin{eqnarray}\label{exe2}
[w_0 \otimes w_1 \otimes L(-1)^nw_2 \otimes w_3] + b_1(z_1, z_2)[w_0 \otimes w_1 \otimes L(-1)^{n-1}w_2 \otimes w_3]\nn
+ \cdots + b_n(z_1, z_2)[w_0 \otimes w_1 \otimes w_2 \otimes w_3] = 0.
\end{eqnarray}
\end{thm}
\pf Without loss of generality, we assume the elements $w_i \in W_i$, $i=0,1,2,3$ are $\N$-weight homogeneous with $\sum_{i=0}^4 \N\text{-}\wt w_i = m$ for some $m \in \N$.

By corollary \ref{maincor}, $F^m(T)/F^m(J)$ is finitely generated. Since $R$ is a Noetherian ring, any $R$-submodule of the finitely generated $R$-module $F^m(T)/F^m(J)$ is also finitely generated. In particular, $M_1$ and $M_2$ are finitely generated. The second conclusion follows immediately. \epfv

Now we establish the existence of systems of differential equations:
\begin{thm}\label{main theorem ex}
Let $W_i$ for $i = 0, 1, 2, 3$ be strongly $\N$-graded generalized $V$-modules satisfying $C_1$-cofiniteness condition with respect to $\N$. Then for any $w_i \in W_i$ ($i = 0, 1, 2, 3$), there exist
\[
a_k(z_1, z_2), b_l(z_1, z_2) \in \C[z_1^{\pm}, z_2^{\pm}, (z_1 - z_2)^{-1}]
\]
for $k = 1, \dots, m$ and $l = 1, \dots, n$ such that for any strongly $\N$-graded $V$-module $W_4, W_5$ and $W_6$, any logarithmic intertwining operators $\mathcal{Y}_1, \mathcal{Y}_2, \mathcal{Y}_3, \mathcal{Y}_4, \mathcal{Y}_5$ and $\mathcal{Y}_6$ of types ${W_0'\choose W_1\,W_4}$, ${W_4\choose W_2\,W_3}$, ${W_5\choose W_1\,W_2}$, ${W_0'\choose W_5\,W_3}$, ${W_0'\choose W_2\,W_6}$ and ${W_6\choose W_1\,W_3}$, respectively, the series
\begin{equation}\label{exe3}
\langle w_0, \mathcal{Y}_1(w_1, z_1)\mathcal{Y}_2(w_2, z_2)w_3\rangle,
\end{equation}
\begin{equation}\label{exe4}
\langle w_0, \mathcal{Y}_4(\mathcal{Y}_3(w_1, z_1-z_2)w_2, z_2)w_3\rangle
\end{equation}
and
\begin{equation}\label{exe5}
\langle w_0, \mathcal{Y}_5(w_2, z_2)\mathcal{Y}_6(w_1, z_1)w_3\rangle,
\end{equation}
satisfy the expansions of the system of differential equations
\begin{equation}\label{exe6}
\frac{\partial^m \varphi}{\partial z_1^m} + a_1(z_1, z_2)\frac{\partial^{m-1} \varphi}{\partial z_1^{m-1}} + \cdots + a_m(z_1, z_2)\varphi = 0,
\end{equation}
\begin{equation}\label{exe7}
\frac{\partial^n \varphi}{\partial z_2^n} + b_1(z_1, z_2)\frac{\partial^{n-1} \varphi}{\partial z_2^{n-1}} + \cdots + b_n(z_1, z_2)\varphi = 0
\end{equation}
in the region $|z_1| > |z_2| > 0$, $|z_2|>|z_1-z_2|>0$ and $|z_2|>|z_1|>0$, respectively.
\end{thm}

\pf The proof is similar to the proof of Theorem 1.4 in \cite{H}. We sketch the proof as follows:

Without loss of generality, we assume the elements $w_i \in W_i$, $i=0,1,2,3$ are $\N$-weight homogeneous with $\sum_{i=0}^4 \N\text{-}\wt w_i = m$ for some $m \in \N$.

Let $\Delta = \wt w_0 - \wt w_1- \wt w_2 - \wt w_3$. Let $\C(\{x\})$ be the space of all series of the form $\sum_{n \in \R}a_nx^n$ for $n \in \R$ such that $a_n = 0$ when the real part of $n$ is sufficiently negative.

Consider the map
\begin{eqnarray*}
\phi_{\mathcal{Y}_1, \mathcal{Y}_2}: F^m(T) \longrightarrow z_1^{\Delta}\C(\{z_2/z_1\})[z_1^{\pm 1}, z_2^{\pm 1}]
\end{eqnarray*}
defined by
\begin{eqnarray*}
& \phi_{\mathcal{Y}_1, \mathcal{Y}_2}(f(z_1, z_2)w_0 \otimes w_1 \otimes w_2 \otimes w_3) \\
& = \iota_{|z_1| > |z_2| > 0}(f(z_1, z_2))\langle w_0, \mathcal{Y}_1(w_1, z_1)\mathcal{Y}_2(w_2, z_2)w_3\rangle,
\end{eqnarray*}
where
\begin{eqnarray*}
\iota_{|z_1| > |z_2| > 0}: R &\longrightarrow & \C[[z_2/z_1]][z_1^{\pm 1}, z_2^{\pm 1}]
\end{eqnarray*}
is the map expanding elements of $R$ as series in the regions $|z_1| > |z_2| > 0$.

Using the Jacobi identity for the logarithmic intertwining operators, we have that $\phi_{\mathcal{Y}_1, \mathcal{Y}_2}(J) = 0$. Thus the map $\phi_{\mathcal{Y}_1, \mathcal{Y}_2}$ induces a map
\begin{eqnarray*}
\bar{\phi}_{\mathcal{Y}_1, \mathcal{Y}_2}: F^m(T)/F^m(J) \longrightarrow z_1^{\Delta}\C(\{z_2/z_1\})[z_1^{\pm 1}, z_2^{\pm 1}].
\end{eqnarray*}
Applying $\bar{\phi}_{\mathcal{Y}_1, \mathcal{Y}_2}$ to (\ref{exe1}) and (\ref{exe2}) and then use the $L(-1)$-derivative property for logarithmic intertwining operators, we see that (\ref{exe3}) indeed satisfies the expansions of the system of differential equations in the regions $|z_1| > |z_2| > 0$. Similarly, we can prove that (\ref{exe4}) and (\ref{exe5}) satisfy the expansions of the system of differential equations in the regions $|z_2| > |z_1-z_2| > 0$ and $|z_2|>|z_1|>0$, respectively. \epfv

\begin{cor}\label{cordiff}Let $W_i = G(\lambda_i, c, l)$ for $|c|<1$ and $\lambda_i \in \mathfrak{h}^{*}$, $i = 0, 1, 2, 3$. Then for any $w_i \in W_i$, there exist
\[
a_k(z_1, z_2), b_l(z_1, z_2) \in \C[z_1^{\pm}, z_2^{\pm}, (z_1 - z_2)^{-1}]
\]
for $k = 1, \dots, m$ and $l = 1, \dots, n$ such that for any strongly $\N$-graded $V$-module $W_4, W_5$ and $W_6$, any logarithmic intertwining operators $\mathcal{Y}_1, \mathcal{Y}_2, \mathcal{Y}_3, \mathcal{Y}_4, \mathcal{Y}_5$ and $\mathcal{Y}_6$ of types ${W_0'\choose W_1\,W_4}$, ${W_4\choose W_2\,W_3}$, ${W_5\choose W_1\,W_2}$, ${W_0'\choose W_5\,W_3}$, ${W_0'\choose W_2\,W_6}$ and ${W_6\choose W_1\,W_3}$, respectively, the series
(\ref{exe3}), (\ref{exe4}), (\ref{exe5}) satisfy the expansions of the system of differential equations
(\ref{exe6}), (\ref{exe7}) in the region $|z_1| > |z_2| > 0$, $|z_2|>|z_1-z_2|>0$ and $|z_2|>|z_1|>0$, respectively.
\end{cor}

Using the fact that the matrix elements of logarithmic intertwining operators among strongly graded modules for $M(l)$ satisfy differential equations (Corollary \ref{cordiff}) and that the singular points for the differential equations are regular (see \cite{Y2}), we have the following convergence and expansion property for $M(l)$ (see \cite{HLZ7}):
\begin{thm}
For any $n \in \Z_{+}$, any strongly graded modules $W_1, \dots, W_{n+1}$ and $\widetilde{W_1}, \dots, \widetilde{W_{n-1}}$, any logarithmic intertwining operators \[\mathcal{Y}_1, \mathcal{Y}_2, \dots, \mathcal{Y}_{n-1}, \mathcal{Y}_n\]of types \[{W_0\choose W_1\,\widetilde{W_1}}, {\widetilde{W_1}\choose W_2\,\widetilde{W_2}}, \dots, {\widetilde{W_{n-2}}\choose W_{n-1}\,\widetilde{W_{n-1}}}, {\widetilde{W_{n-1}}\choose W_n\,W_{n+1}},\]respectively, and any $w_{(0)}' \in W_0'$, $w_{(1)} \in W_1, \dots, W_{(n+1)} \in W_{n+1}$, the series \[\langle w_{(0)}', \mathcal{Y}_1(w_{(1)}, z_1)\cdots \mathcal{Y}_n(w_{(n)}, z_n)w_{(n+1)}\rangle\]is absolutely convergent in the region $|z_1| > \cdots > |z_n| > 0$ and its sum can be analytically extended to a multivalued analytic function on the region given by $z_1 \neq 0$, $i = 1, \dots, n$, $z_i \neq z_j$, $i \neq j$, such that for any set of possible singular points with either $z_i = 0, z_i = \infty$ or $z_i = z_j$ for $i \neq j$, this multivalued analytic function can be expanded near the singularity as a series having the same form as the expansion near the singular points of a solution of a system of differential equations with regular singular points.
\end{thm}
\pf The proof uses the same method as in \cite{H}. \epfv

\section{Graded dimension}
Let $V = \coprod_{\alpha \in A, n \in \Z}V_{(n)}^{(\alpha)}$ be a strongly $A$-graded vertex algebra. We call the {\it graded dimension} of $V$ the formal sum
\[
\mbox{dim}_{*}\;V = \sum_{\alpha \in A, n \in \Z}(\mbox{dim}\; V_{(n)}^{(\alpha)})p^{n}q^{\alpha}.
\]
Here $p$ and $q$ are mutually commuting formal variables. In particular, the graded dimension of $M(l)$ is
\[
\mbox{dim}_{*}\;M(l) = \sum_{m, n \in \N}\mbox{dim}\;(M(l)_{(n)}^{(m)})p^nq^m.
\]

The dimension of the double homogeneous subspace $M(l)_{(n)}^{(m)}$ is the number of bipartite partitions
\[
(m, n) = \sum_{j = 1}^{\infty}(m_j, n_j),
\]
where $m_j$'s are nonnegative integers and $n_j$'s are positive integers for $j \in \Z_{+}$. We use the following standard abbreviation in \cite{A}:
\[
(a; q)_n = (1-a)(1-aq)\cdots (1-aq^{n-1}),
\]
\[
(a; q)_{\infty} = \lim_{n \mapsto \infty}(a; q)_n,
\]
where $a$ is any formal variable or complex number.

Let $p(k, n)$ be the number of partitions of $n$ into a set of positive integers that have $k$ parts. Then (see for example \cite{A})
\[
\sum_{k = 0}^{\infty}\sum_{n = 0}^{\infty}p(k,n)x^kq^n = \frac{1}{(xq; q)_{\infty}}.
\]
Similarly, let $p'(k, m)$ the number of partitions of $m$ into a set of nonnegative integers that have $k$ parts. Then
\[
\sum_{k = 0}^{\infty}\sum_{m = 0}^{\infty}p'(k,m)x^kq^m = \frac{1}{(x; q)_{\infty}}.
\]
Thus we have
\begin{equation}
\mbox{dim}_{*}\;M(l) = \mbox{constant term of $x$ in}\; \frac{1}{(x^{-1}p; p)_{\infty}(x; q)_{\infty}}.
\end{equation}

\def\refname{\hfil{REFERENCES}}

\noindent {\small \sc Department of Mathematics, University of Notre Dame,
255 Hurley Building, Notre Dame, IN 46556-4618}
\vspace{1em}

\noindent {\em E-mail address}: jyang7@nd.edu

\end{document}